# Fréchet and Mordukhovich Derivative (Coderivative) and Covering Constant for Single-Valued Mapping in Euclidean Space with Applications (I)

Jinlu Li

Department of Mathematics
Shawnee State University
Portsmouth, OH 45662 USA
Email: jli@shawnee.edu

**Abstract** In this paper, we study Fréchet derivatives and Mordukhovich derivatives (or coderivatives) of single-valued mappings in Euclidean spaces. At first, we prove the guideline for calculating the Fréchet derivatives of single-valued mappings by their partial derivatives. Then, by using the connections between Fréchet derivatives and Mordukhovich derivatives (or coderivatives) of single-valued mappings in Banach spaces, we derive the useful rules for calculating the Mordukhovich derivatives of single-valued mappings in Euclidean spaces. For practicing these rules, we find the precise solutions of the Fréchet derivatives and Mordukhovich derivatives for some single-valued mappings in Euclidean spaces (in $\mathbb{R}^2$, it can be extended to $\mathbb{R}^n$). By using these solutions, we will find the covering constants for the considered mappings. As applications of the results about the covering constants and by applying the Arutyunov Mordukhovich and Zhukovskiy Parameterized Coincidence Point Theorem, we solve some parameterized equations.

## 1. Introduction

The concept of Fréchet derivative of a single-valued mapping is a generalization of the derivative of a real-valued function of a single real variable in Euclidean spaces in calculus to the case of a vector-valued function of multiple real variables in normed spaces, in particular, in Banach spaces. The Fréchet derivative has been widely used in the calculus of variations and nonlinear functional analysis in Banach spaces (See [8, 9, 11, 22]).

One more step further, in the theory of generalized differentiational analysis in Banach spaces, the concept of Mordukhovich derivative (or coderivative) for set-valued mappings lays the foundation and it plays fundamental and crucial role in the theory of set-valued and variational analysis (See [18−21, 23]). The results of Theorem 1.38 in [18] provides close connections between the Fréchet derivative and the Mordukhovich derivative (or coderivative) of single-valued mappings in Banach spaces. By this theorem, the Mordukhovich derivative can be considered to be the generalization of the Fréchet derivative of single-valued mappings.

The theory of generalized differentiation in set-valued analysis has been rapidly developed during the last years and has been widely applied to many fields such as optimization, control theory, game theory, variational analysis, and so forth (see [1−7, 12−21, 24]). It is well-known that, in contrast to the differentiation theory in calculus when the generalized differentiation is applied to a specific problem with a specific mapping, one needs to find the exact solution of the Mordukhovich derivative of the considered mapping. For example, the concept of covering constant plays an important and crucial role in Arutyunov Mordukhovich and Zhukovskiy Parameterized Coincidence Point Theorem. However, the covering constants are defined by the Mordukhovich derivatives. So, if one wants to find the covering constant for a mapping, you have firstly to find the Mordukhovich derivatives of the considered mapping. This idea leads an important motivation for authors to find some formulas for both the Fréchet derivatives and Mordukhovich derivatives of mappings in Banach spaces. For example, in [12−13], the precise solutions for the standard metric project operator have been proved in both Hilbert spaces and uniformly convex and uniformly smooth Banach spaces.

In this paper, we concentrate to study the Fréchet derivatives and Mordukhovich derivatives of single-valued mappings in Euclidean spaces, which is implemented as the following steps.

(i) At first, we prove the principle of the Fréchet differentiability of single-valued mappings in Euclidean spaces.
(ii) By the second order approximations of mappings, we derive the guidelines for calculating the Fréchet derivatives.
(iii) We use the relationship between Fréchet derivatives and Mordukhovich derivatives (It is given by Theorem 1.38 in [18]) to deduce the algorithm for calculating the Mordukhovich derivatives of single-valued mappings in Euclidean spaces.
(iv) We provide examples of mappings in $\mathbb{R}^2$, which includes polynomial mappings, rational mappings, exponential mappings and logarithm mappings. For these mappings, we find the

   (a) Fréchet derivatives,
   (b) Mordukhovich derivatives,
   (c) Covering constants.

(v) By the results of covering constants obtained in Step 4, we apply the Arutyunov Mordukhovich and Zhukovskiy Parameterized Coincidence Point Theorem to solve some parameterized equations.

2. **Review for Fréchet derivative and Mordukhovich derivative (coderivative)**

In this section, we briefly review the concepts and properties of Fréchet derivatives and Mordukhovich derivatives (or coderivatives) of single-valued mappings in Banach spaces, in particular, in Euclidean spaces, which will be used in this paper. See [18−21] for more details,

Let $(X, \|\cdot\|_X)$ and $(Y, \|\cdot\|_Y)$ be real Banach spaces, which have topological dual spaces $(X^*, \|\cdot\|_{X^*})$ and $(Y^*, \|\cdot\|_{Y^*})$, respectively. Let $\langle \cdot, \cdot \rangle_X$ denote the real canonical pairing between $X^*$ and $X$ and $\langle \cdot, \cdot \rangle_Y$ the real canonical pairing between $Y^*$ and $Y$. Let $\theta_X$ and $\theta_Y$ denote the origins in $X$ and $Y$, respectively. For any $x \in X$ and $r > 0$, let $B_X(x, r)$ and $S_X(x, r)$ denote the closed ball and sphere in $X$ centered at point $x$ with radius $r$, respectively.

Let $g: X \to Y$ be a single-valued mapping and let $\bar{x} \in X$. If there is a linear and continuous mapping $\nabla g(\bar{x}): X \to Y$ such that

$$\lim_{x \to \bar{x}} \frac{g(x) - g(\bar{x}) - \nabla g(\bar{x})(x - \bar{x})}{\|x - \bar{x}\|_X} = \theta_Y,$$

then $g$ is said to be Fréchet differentiable at $\bar{x}$ and $\nabla g(\bar{x})$ is called the Fréchet derivative of $g$ at $\bar{x}$.

The Mordukhovich derivative (coderivative) is initially defined for set-valued mappings in Banach spaces. Since in this paper, we only deal with the single-valued mappings in Euclidean spaces, so we only review the concept of Mordukhovich derivative (coderivative) for single-valued mappings in Banach spaces (See [18–21] for more details). Let $\Delta$ be a nonempty subset of $X$ and let $g: \Delta \to Y$ be a single-valued mapping. For $x \in \Delta$ and $y = g(x)$, a set-valued mapping $\widehat{D}^*g(x,y): Y^* \to X^*$ is defined, for any $y^* \in Y^*$, by (See Definitions 1.13 and 1.32 in Chapter 1 in [15])

$$\widehat{D}^*g(x,y)(y^*) = \left\{z^* \in X^*: \limsup_{\substack{(u,g(u)) \to (x,g(x)) \\ u \in \Delta}} \frac{\langle z^*, u-x \rangle_X - \langle y^*, g(u)-g(x) \rangle_Y}{\|u-x\|_X + \|g(u)-g(x)\|_Y} \leq 0\right\}.$$

If this set-valued mapping $\widehat{D}^*g(x,y): Y^* \to X^*$ satisfies that

$$\widehat{D}^*g(x,y)(y^*) \neq \emptyset, \text{ for every } y^* \in Y^*, \tag{2.1}$$

then $g$ is said to be Mordukhovich differentiable (codifferentiable) at the point $x$ and $\widehat{D}^*g(x,y)$ is called the Mordukhovich derivative (which is also called Mordukhovich coderivative, or coderivative) of $g$ at point $x$. For this single-valued mapping $g: \Delta \to Y$, we write

$$\widehat{D}^*g(x,y)(y^*) \equiv \widehat{D}^*g(x)(y^*).$$

Furthermore, if $g: \Delta \to Y$ is a single-valued continuous mapping. Then, by the above definition, the Mordukhovich derivative of $g$ at point $x$ is calculated by

$$\widehat{D}^*g(x)(y^*) = \left\{z^* \in X^*: \limsup_{\substack{u \to x \\ u \in \Delta}} \frac{\langle z^*, u-x \rangle_X - \langle y^*, g(u)-g(x) \rangle_Y}{\|u-x\|_X + \|g(u)-g(x)\|_Y} \leq 0\right\}, \text{ for any } y^* \in Y^*. \tag{2.2}$$

The following theorem shows the connection between Fréchet derivatives and Mordukhovich derivatives for sing-valued mappings. The results of the following theorem provide a powerful tool to calculate the Mordukhovich derivatives by the Fréchet derivatives of single-valued mappings.

**Theorem 1.38 in [15]**. *Let $X$ be a Banach space with dual space $X^*$ and let $g: X \to Y$ be a single-valued mapping. Suppose that $g$ is Fréchet differentiable at $x \in X$ with $y = g(x)$. Then, the Mordukhovich derivative of $g$ at $x$ satisfies the following equation*

$$\widehat{D}^*g(x)(y^*) = \{(\nabla g(x))^*(y^*)\}, \text{ for all } y^* \in Y^*.$$

One of the important applications of Mordukhovich derivatives of set-valued mappings is to define the covering constants for set-valued mappings. The covering constant for $\Phi: X \rightrightarrows Y$ at point $(\bar{x}, \bar{y}) \in \text{gph } \Phi$ is defined by (see (2.4) in [1])

$$\hat{\alpha}(\Phi, \bar{x}, \bar{y}) := \sup_{\eta > 0} \inf\{\|z^*\|_{X^*}: z^* \in \widehat{D}^*\Phi(x,y)(w^*), x \in \mathbb{B}_X(\bar{x}, \eta), y \in \Phi(x) \cap \mathbb{B}_Y(\bar{y}, \eta), \|w^*\|_{Y^*} = 1\}. \tag{2.4}$$

Here, $\|\cdot\|_{X^*}$ and $\|\cdot\|_{Y^*}$ denote the norms in $X^*$ and $Y^*$, respectively. $\mathbb{B}_X(\bar{x}, \eta)$ is the closed ball in $X$ centered at $\bar{x}$ with radius $\eta$, and $\mathbb{B}_Y(\bar{y}, \eta)$ is the closed ball in $Y$ centered at $\bar{y}$ with radius $\eta$.

In particular, let $g: X \to Y$ be a single-valued mapping. For any $\bar{x}, \bar{y} \in X$ with $\bar{y} = g(\bar{x})$, (2.4) becomes

$$\hat{\alpha}(g, \bar{x}, \bar{y}) = \sup_{\eta > 0} \inf\{\|z^*\|_{X^*}: z^* \in \widehat{D}^*g(x)(w^*), x \in \mathbb{B}_X(\bar{x}, \eta), g(x) \in \mathbb{B}_Y(\bar{y}, \eta), \|w^*\|_{Y^*} = 1\}. \tag{2.5}$$

In particular, in this paper, we consider Hilbert spaces as special cases of Banach spaces. Let $(H, \|\cdot\|)$ be a real Hilbert space with inner product $\langle \cdot, \cdot \rangle$ and origin $\theta$. For $r > 0$ and $x \in H$, let $\mathbb{B}(x, r)$ be the closed ball in $H$ with radius $r$ and centered at $x$. Let $g: H \to H$ be a single-valued mapping and let $\bar{x} \in X$. Then, the Fréchet derivative of $g$ at $\bar{x}$ is defined by $\nabla g(\bar{x}): H \to H$ such that

$$\lim_{x \to \bar{x}} \frac{g(x) - g(\bar{x}) - \nabla g(\bar{x})(x - \bar{x})}{\|x - \bar{x}\|} = \theta.$$

Let $\Delta$ be a nonempty subset of $H$ and let $g: \Delta \to H$ be a single-valued continuous mapping. For $x \in \Delta$ and $y = g(x)$, by (2.2), the Mordukhovich derivative of $g$ at point $x$ becomes

$$\widehat{D}^* g(x)(y) = \left\{ z \in H : \limsup_{\substack{u \to x \\ u \in \Delta}} \frac{\langle z, u - x \rangle - \langle y, g(u) - g(x) \rangle}{\|u - x\| + \|g(u) - g(x)\|} \leq 0 \right\}, \text{ for any } y \in H. \tag{2.6}$$

In this case, for any $\bar{x}, \bar{y} \in H$ with $\bar{y} = g(\bar{x})$, by (2.4), the covering constant for $g$ at point $(\bar{x}, \bar{y})$

$$\hat{a}(g, \bar{x}, \bar{y}) = \sup_{\eta > 0} \inf \{\|z\|: z \in \widehat{D}^* g(x)(w), x \in \mathbb{B}(\bar{x}, \eta), f(x) \in \mathbb{B}(\bar{y}, \eta), \|w\| = 1\}. \tag{2.7}$$

### 3. Fréchet and Mordukhovich Derivatives (Coderivatives) of Mappings in $\mathbb{R}^n$

Let $n \geq 1$ and let $(\mathbb{R}^n, \|\cdot\|)$ be the standard $n$-d Euclidean space with the ordinal Hilbert $L_2$-norm and row vectors. Let $\theta$ denote the origin of $\mathbb{R}^n$. Let $m, n \geq 1$ and let $f: \mathbb{R}^n \to \mathbb{R}^m$ be a single-valued mapping with the following representation.

$$f((x_1, x_2, \ldots, x_n)) = \big(f_1(x_1, x_2, \ldots, x_n), f_2(x_1, x_2, \ldots, x_n), \ldots, f_m(x_1, x_2, \ldots, x_n)\big), \text{ for } (x_1, x_2, \ldots, x_n) \in \mathbb{R}^n.$$

Where, for $i = 1, 2, \ldots, m$, $f_i(x_1, x_2, \ldots, x_n): \mathbb{R}^n \to \mathbb{R}$ is a real valued multivariable function defined on $\mathbb{R}^n$ with variables $x_1, x_2, \ldots, x_n$. We write the above equation as $f = (f_1, f_2, \ldots, f_m)$. Let $(z_1, z_2, \ldots, z_n) \in \mathbb{R}^n$. For $i = 1, 2, \ldots, n$, the partial derivative of $f_i$ with respect to the variable $x_j$ at $(z_1, z_2, \ldots, z_n)$ is

$$\frac{\partial f_i}{\partial x_j}(z_1, z_2, \ldots, z_n), \text{ for every } j = 1, 2, \ldots, n. \tag{3.1}$$

By equation (3.1), for any $i = 1, 2, \ldots, m$, the existence of $f_i(z_1, z_2, \ldots, z_n)$ means that,

$$\lim_{x_j \to z_j} \frac{f_i(x_1, x_2, \ldots, x_n) - f_i(z_1, z_2, \ldots, z_n) - \frac{\partial f_i}{\partial x_j}(z_1, z_2, \ldots, z_n)(x_j - z_j)}{x_j - z_j} = 0, \text{ for } j = 1, 2, \ldots, n. \tag{3.2}$$

The limit (3.2) is equivalent to, for any $i = 1, 2, \ldots, m$,

$$\lim_{x_j \to z_j} \frac{f_i(x_1, x_2, \ldots, x_n) - f_i(z_1, z_2, \ldots, z_n)}{x_j - z_j} = -\frac{\partial f_i}{\partial x_j}(z_1, z_2, \ldots, z_n), \text{ for } j = 1, 2, \ldots, n. \tag{3.3}$$

Furthermore, if $f$ satisfies the following conditions

$$\lim_{x \to z} \frac{f_i(x_1, x_2, \ldots, x_n) - f_i(z_1, z_2, \ldots, z_n) - \left(\sum_{j=1}^n \frac{\partial f_i}{\partial x_j}(z_1, z_2, \ldots, z_n)(x_j - z_j)\right)}{\|x - z\|} = 0, \text{ for } i = 1, 2, \ldots, m, \tag{3.4}$$

then, $f$ is differentiable and has the linear approximation at point $z = (z_1, z_2, \ldots, z_n)$, which is the first order of the Taylor polynomial of the mapping $f$ from $\mathbb{R}^n$ to $\mathbb{R}^m$. We have some sufficient conditions for $f$ to have the linear approximation at point $z$.

For every $i = 1, 2, \ldots, m$, the second order partial derivative of $f_i$ at point $z = (z_1, z_2, \ldots, z_n)$ with respect to $x_j$ and $x_k$ is

$$\frac{\partial^2 f_i}{\partial x_j \partial x_k}(z_1, z_2, \ldots, z_n), \text{ for } j, k = 1, 2, \ldots, n.$$

**Fact 3.1.** *Let $f = (f_1, f_2, \ldots, f_m). : \mathbb{R}^n \to \mathbb{R}^m$ be a single-valued mapping. Let $z \in \mathbb{R}^n$. If there is a ball $B$ with radius $r > 0$ and centered at $z$ such that, for every $i = 1, 2, \ldots, m$, the real valued function $f_i: \mathbb{R}^n \to \mathbb{R}$ is twice differentiable (all second partial derivatives of $f_i$ exist) at every point $y \in B$, that is*

$$\frac{\partial^2 f_i}{\partial x_j \partial x_k}(z_1, z_2, \ldots, z_n) \text{ exists for any } z \in B, \text{ for } j, k = 1, 2, \ldots, n.$$

*then $f$ has the linear approximation at this point $z$.*

*Proof.* This is a well-known result. It is proved by using the second order Taylor polynomial of each real valued function $f_i: \mathbb{R}^n \to \mathbb{R}$, for $i = 1, 2, \ldots, n$. The details are omitted here. □

**Theorem 3.2.** *Let $f = (f_1, f_2, \ldots, f_m): \mathbb{R}^n \to \mathbb{R}^m$ be a single-valued mapping. Let $(z_1, z_2, \ldots, z_n) \in \mathbb{R}^n$. Suppose that, for every $i = 1, 2, \ldots, m$, $\frac{\partial f_i}{\partial x_j}(z_1, z_2, \ldots, z_n)$ exists, for every $j = 1, 2, \ldots, n$ and $f$ has the linear approximation at point $z$. Then,*

(a) *$f$ is Fréchet differentiable at $z$ and the Fréchet derivative of $f$ at $z$ is the following $n \times m$ matrix,*

$$\nabla f(z) = \begin{pmatrix} \frac{\partial f_1}{\partial x_1}(z_1, z_2, \ldots, z_n) & \cdots & \frac{\partial f_m}{\partial x_1}(z_1, z_2, \ldots, z_n) \\ \vdots & \ddots & \vdots \\ \frac{\partial f_1}{\partial x_n}(z_1, z_2, \ldots, z_n) & \cdots & \frac{\partial f_m}{\partial x_n}(z_1, z_2, \ldots, z_n) \end{pmatrix}, \quad (3.5)$$

(b) *$f$ is Mordukhovich differentiable at point $z$ that is the Jacobian matrix of $f$ at $z$*

$$\widehat{D}^* f(x) = \nabla f(z)^T = \begin{pmatrix} \frac{\partial f_1}{\partial x_1}(z_1, z_2, \ldots, z_n) & \cdots & \frac{\partial f_1}{\partial x_n}(z_1, z_2, \ldots, z_n) \\ \vdots & \ddots & \vdots \\ \frac{\partial f_m}{\partial x_1}(z_1, z_2, \ldots, z_n) & \cdots & \frac{\partial f_m}{\partial x_n}(z_1, z_2, \ldots, z_n) \end{pmatrix}. \quad (3.6)$$

*Proof.* Proof of (a). Notice that, in this paper, $\mathbb{R}^n$ has row vectors. By the meanings of (3.2) and (3.3), if $f$ has the linear approximation at this point $z$, by (3.4), we calculate

$$\lim_{x \to z} \frac{f(x_1, x_2, \ldots, x_n) - f(z_1, z_2, \ldots, z_n) - \nabla f(z)(x-z)}{\|x-z\|}$$

$$= \lim_{x \to z} \left( \frac{(f_1(x_1, x_2, \ldots, x_n), f_2(x_1, x_2, \ldots, x_n), \ldots, f_m(x_1, x_2, \ldots, x_n)) - (f_1(z_1, z_2, \ldots, z_n), f_2(z_1, z_2, \ldots, z_n), \ldots, f_m(z_1, z_2, \ldots, z_n))}{\|x-z\|} \right.$$

$$-\frac{(x_1-z_1,\ldots,x_n-z_n)\begin{pmatrix}\frac{\partial f_1}{\partial x_1}(z_1,z_2,\ldots,z_n) & \cdots & \frac{\partial f_m}{\partial x_1}(z_1,z_2,\ldots,z_n)\\ \vdots & \ddots & \vdots\\ \frac{\partial f_1}{\partial x_n}(z_1,z_2,\ldots,z_n) & \cdots & \frac{\partial f_m}{\partial x_n}(z_1,z_2,\ldots,z_n)\end{pmatrix}}{\|x-z\|}$$

$$=\lim_{x\to z}\left(\frac{f_1(x_1,x_2,\ldots,x_n)-f_1(z_1,z_2,\ldots,z_n)-\left(\sum_{j=1}^n (x_j-z_j)\frac{\partial f_1}{\partial x_j}(z_1,z_2,\ldots,z_n)\right)}{\|x-z\|},\ldots,\right.$$

$$\left.\frac{f_m(x_1,x_2,\ldots,x_n)-f_m(z_1,z_2,\ldots,z_n)-\left(\sum_{j=1}^n (x_j-z_j)\frac{\partial f_m}{\partial x_1}(z_1,z_2,\ldots,z_n)\right)}{\|x-z\|}\right)$$

$= (0, \ldots, 0)$.

Part (b) can be proved by Theorem 1.38 in [15] and Part (a) of this theorem. □

**Proposition 3.3.** *Let $f = (f_1, f_2, \ldots, f_m): \mathbb{R}^n \to \mathbb{R}^m$ be a single-valued mapping.*

(a) *Suppose that, for every $i = 1, 2, \ldots, m$, $f_i(x_1, x_2, \ldots, x_n)$ is a polynomial function with respect to $x_1, x_2, \ldots, x_n$, then $f$ is Fréchet differentiable at every point in $\mathbb{R}^n$, and $f$ is Mordukhovich differentiable on $\mathbb{R}^n$;*

(b) *Suppose that, for every $i = 1, 2, \ldots, m$, $f_i(x_1, x_2, \ldots, x_n)$ is a rational function with respect to $x_1, x_2, \ldots, x_n$. Let $z \in \mathbb{R}^n$. If $z$ is not a zero point for the denominator of every $f_i$, then $f$ is Fréchet differentiable at $z$; and therefore, $f$ is Mordukhovich differentiable at $z$.*

## 4. Calculating Covering Constants for Single-Valued Mappings in $\mathbb{R}^2$

In this section, as applications of Lemma 3.1 and Theorems 3.2 and 3.3 in the previous section, we provide some examples to calculate the Fréchet and Mordukhovich derivatives and covering constants of single-valued mappings in $\mathbb{R}^2$, which can be extended to higher dimension Euclidean spaces. In the examples in this section, we will show the details and techniques for calculating the covering constants for single-valued mappings in $\mathbb{R}^2$. From these examples, we can see the difficulty and complexity for the calculation of covering constants.

**A rational mapping 4.1.** We define $f \equiv (f_1, f_2): \mathbb{R}^2 \to \mathbb{R}^2$, for $x = (x_1, x_2) \in \mathbb{R}^2$, by

$$f(x) = \begin{cases}\left(\frac{x_1^2-x_2^2}{x_1^2+x_2^2}, \frac{2x_1x_2}{x_1^2+x_2^2}\right), & \text{for } (x_1, x_2) \neq \theta,\\ \theta, & \text{for } (x_1, x_2) = \theta.\end{cases}$$

Then, $f$ has the following properties.

(a) *$f$ is Fréchet differentiable and Mordukhovich differentiable on $\mathbb{R}^2\setminus\{\theta\}$. For any $z = (z_1, z_2) \in \mathbb{R}^2\setminus\{\theta\}$, we have*

$$\nabla f(z) = \begin{pmatrix} \frac{4z_1 z_2^2}{(z_1^2+z_2^2)^2} & \frac{-2z_2(z_1^2-z_2^2)}{(z_1^2+z_2^2)^2} \\ \frac{-4z_1^2 z_2}{(z_1^2+z_2^2)^2} & \frac{2z_1(z_1^2-z_2^2)}{(z_1^2+z_2^2)^2} \end{pmatrix},$$

and
$$\widehat{D}^* f(z) = \begin{pmatrix} \frac{4z_1 z_2^2}{(z_1^2+z_2^2)^2} & \frac{-4z_1^2 z_2}{(z_1^2+z_2^2)^2} \\ \frac{-2z_2(z_1^2-z_2^2)}{(z_1^2+z_2^2)^2} & \frac{2z_1(z_1^2-z_2^2)}{(z_1^2+z_2^2)^2} \end{pmatrix},$$

(b) The covering constant for $f$ is constant on $\mathbb{R}^2\backslash\{\theta\}$ satisfying

$$\hat{\alpha}(f,\bar{z},\bar{w}) = 0, \text{ for any } \bar{z} = (\bar{z}_1,\bar{z}_2) \neq \theta \text{ with } \bar{w} = f(\bar{z}).$$

(c) $f$ is not Fréchet differentiable at $\theta$.

(d) $f$ is not Mordukhovich differentiable at $\theta$. More precisely speaking, we have

$$\widehat{D}^* f(\theta)(y) = \emptyset, \text{ for any } y \neq \theta.$$

*Proof.* Proof of (a). This mapping $f$ a rational mapping and it is continuous on $\mathbb{R}^2\backslash\{\theta\}$. We calculate the norms of this mapping.

$$\|f(x_1,x_2)\|^2 = \left(\frac{x_1^2-x_2^2}{x_1^2+x_2^2}\right)^2 + \left(\frac{2x_1 x_2}{x_1^2+x_2^2}\right)^2 = 1, \text{ for any } (x_1,x_2) \in \mathbb{R}^2\backslash\{\theta\}$$

This implies that the norm of $f \equiv (f_1, f_2)$ is constant on $\mathbb{R}^2\backslash\{\theta\}$ with

$$\|f(x)\| = 1, \text{ for any } x \in \mathbb{R}^2\backslash\{\theta\}.$$

By Lemma 3.1 and Theorem 3.2, $f$ is Fréchet differentiable and Mordukhovich differentiable at every point $z = (z_1, z_2) \in \mathbb{R}^2\backslash\{\theta\}$. We have

$$\frac{\partial\left(\frac{x_1^2-x_2^2}{x_1^2+x_2^2}\right)}{\partial z_1} = \frac{4z_1 z_2^2}{(z_1^2+z_2^2)^2}, \quad \frac{\partial\left(\frac{2x_1 x_2}{x_1^2+x_2^2}\right)}{\partial z_1} = \frac{-2z_2(z_1^2-z_2^2)}{(z_1^2+z_2^2)^2},$$

$$\frac{\partial\left(\frac{x_1^2-x_2^2}{x_1^2+x_2^2}\right)}{\partial z_2} = \frac{-4z_1^2 z_2}{(z_1^2+z_2^2)^2}, \quad \frac{\partial\left(\frac{2x_1 x_2}{x_1^2+x_2^2}\right)}{\partial z_2} = \frac{2z_1(z_1^2-z_2^2)}{(z_1^2+z_2^2)^2}.$$

This proves (a). Then, we prove (b).

Proof of (b). Let $x = (x_1, x_2)$ and $y = (y_1, y_2) \in \mathbb{R}^2$. If $x = \widehat{D}^* f(z)(y)$, by part (a), we have that

$$(x_1, x_2) = (y_1, y_2) \begin{pmatrix} \frac{4z_1 z_2^2}{(z_1^2+z_2^2)^2} & \frac{-4z_1^2 z_2}{(z_1^2+z_2^2)^2} \\ \frac{-2z_2(z_1^2-z_2^2)}{(z_1^2+z_2^2)^2} & \frac{2z_1(z_1^2-z_2^2)}{(z_1^2+z_2^2)^2} \end{pmatrix}.$$

This is rewritten, for $x = \widehat{D}^* f(z)(y)$, by

$$x_1 = y_1 \frac{4z_1 z_2^2}{(z_1^2+z_2^2)^2} - y_2 \frac{2z_2(z_1^2-z_2^2)}{(z_1^2+z_2^2)^2},$$

$$x_2 = -y_1 \frac{4z_1^2 z_2}{(z_1^2+z_2^2)^2} + y_2 \frac{2z_1(z_1^2-z_2^2)}{(z_1^2+z_2^2)^2}.$$

Hence, if $x = \widehat{D}^* f(z)(y)$, then

$$x_1^2 + x_2^2 = \left(y_1 \frac{4z_1 z_2^2}{(z_1^2+z_2^2)^2} - y_2 \frac{2z_2(z_1^2-z_2^2)}{(z_1^2+z_2^2)^2}\right)^2 + \left(-y_1 \frac{4z_1^2 z_2}{(z_1^2+z_2^2)^2} + y_2 \frac{2z_1(z_1^2-z_2^2)}{(z_1^2+z_2^2)^2}\right)^2$$

$$= y_1^2 \left(\frac{4z_1 z_2^2}{(z_1^2+z_2^2)^2}\right)^2 - 2y_1 y_2 \frac{4z_1 z_2^2}{(z_1^2+z_2^2)^2} \frac{2z_2(z_1^2-z_2^2)}{(z_1^2+z_2^2)^2} + y_2^2 \left(\frac{2z_2(z_1^2-z_2^2)}{(z_1^2+z_2^2)^2}\right)^2$$

$$+ y_1^2 \left(\frac{4z_1^2 z_2}{(z_1^2+z_2^2)^2}\right)^2 - 2y_1 y_2 \frac{4z_1^2 z_2}{(z_1^2+z_2^2)^2} \frac{2z_1(z_1^2-z_2^2)}{(z_1^2+z_2^2)^2} + y_2^2 \left(\frac{2z_1(z_1^2-z_2^2)}{(z_1^2+z_2^2)^2}\right)^2$$

$$= \frac{4}{(z_1^2+z_2^2)^4} \left(4y_1^2 z_1^2 z_2^2 (z_1^2+z_2^2) - 4y_1 y_2 (z_1 z_2 (z_1^2-z_2^2)(z_1^2+z_2^2)) + y_2^2((z_1^2-z_2^2)^2(z_1^2+z_2^2))\right)$$

$$= \frac{4}{(z_1^2+z_2^2)^3} \left(4y_1^2 z_1^2 z_2^2 - 4y_1 y_2 (z_1 z_2 (z_1^2-z_2^2)) + y_2^2 (z_1^2-z_2^2)^2\right)$$

$$= \frac{4(2y_1 z_1 z_2 - y_2(z_1^2-z_2^2))^2}{(z_1^2+z_2^2)^3}.$$

This implies

$$\|x\| = 2 \frac{|2y_1 z_1 z_2 - y_2(z_1^2-z_2^2)|}{\|z\|^3}, \text{ for } x = \widehat{D}^* f(z)(y). \tag{4.1}$$

For any given point $\bar{z} = (\bar{z}_1, \bar{z}_2) \neq \theta$ and $\bar{w} = (\bar{w}_1, \bar{w}_2) \in \mathbb{R}^2$ with $\bar{w} = f(\bar{z})$, by (4.1), we calculate

$$\hat{\alpha}(f, \bar{z}, \bar{w}) = \sup_{\eta > 0} \inf\{\|x\|: x \in \widehat{D}^* f(z)(y), z \in \mathbb{B}(\bar{z}, \eta), f(z) \in \mathbb{B}(\bar{w}, \eta), \|y\| = 1\}.$$

By the continuity of $f$ on $\mathbb{R}^2 \setminus \{\theta\}$, for any given $\eta > 0$, we have

$$\{z \in \mathbb{B}(\bar{z}, \eta), f(z) \in \mathbb{B}(\bar{w}, \eta)\} \neq \emptyset.$$

For any given $\eta > 0$, and $z = (z_1, z_2) \in \mathbb{B}(\bar{z}, \eta) \setminus \{\theta\}$ and $f(z) \in \mathbb{B}(\bar{w}, \eta)$, the following system of equations has solutions with respect to $y = (y_1, y_2)$:

$$\begin{cases} y_1^2 + y_2^2 = 1, \\ 2y_1 z_1 z_2 - y_2(z_1^2 - z_2^2) = 0. \end{cases} \tag{4.2}$$

By (4.1) and (4.2), we have

$$\hat{\alpha}(f, \bar{z}, \bar{w}) = \sup_{\eta > 0} \inf\{\|x\|: x \in \widehat{D}^* f(z)(y), z \in \mathbb{B}(\bar{z}, \eta), f(z) \in \mathbb{B}(\bar{w}, \eta), \|y\| = 1\}$$

$$\leq \sup_{\eta > 0} \inf\{\|x\|: x \in \widehat{D}^* f(z)(y), z \in \mathbb{B}(\bar{z}, \eta) \setminus \{\theta\}, f(z) \in \mathbb{B}(\bar{w}, \eta), \|y\| = 1, 2y_1 z_1 z_2 - y_2(z_1^2 - z_2^2) = 0\}$$

$$= \sup_{\eta>0} \inf\{0: x \in \widehat{D}^*f(z)(y), z \in \mathbb{B}(\bar{z},\eta)\setminus\{\theta\}, f(z) \in \mathbb{B}(\bar{w},\eta), \|y\| = 1, 2y_1 z_1 z_2 - y_2(z_1^2 - z_2^2) = 0\}$$

$$= 0.$$

This proves that

$$\hat{\alpha}(f,\bar{z},\bar{w}) = 0, \text{ for any } \bar{z} = (\bar{z}_1,\bar{z}_2) \neq \theta \text{ with } \bar{w} = f(\bar{z}).$$

Part (b) is proved. It is clear that $f$ is not continuous at $\theta$, which implies that $f$ is not Fréchet differentiable at $\theta$. This shows (c).

Proof of (d). Part (d) can also be proved by the discontinuity of $f$ at $\theta$. Since the Mordukhovich derivative of $f$ is defined by an inequality, we give the proof of (d) here.

Let $y = (y_1, y_2) \in \mathbb{R}^2\setminus\{\theta\}$ and $x = (x_1, x_2) \in \mathbb{R}^2$, we calculate the following limits.

$$\limsup_{u\to\theta} \frac{\langle (x_1,x_2),\ (u_1,u_2)-\theta\rangle - \langle (y_1,y_2),\ f(u)-f(\theta)\rangle}{\|u-\theta\| + \|f(u)-f(\theta)\|}$$

$$= \limsup_{u\to\theta} \frac{\langle (x_1,x_2),\ (u_1,u_2)\rangle - \langle (y_1,y_2),\ \left(\frac{u_1^2-u_2^2}{u_1^2+u_2^2}, \frac{2u_1 u_2}{u_1^2+u_2^2}\right)\rangle}{\|u\| + \|f(u)\|}$$

$$= \limsup_{u\to\theta} \frac{x_1 u_1 + x_2 u_2 - y_1 \frac{u_1^2-u_2^2}{u_1^2+u_2^2} - y_2 \frac{2u_1 u_2}{u_1^2+u_2^2}}{\|u\|+1}. \tag{4.3}$$

Suppose $y_1 \neq 0$. We consider the following four cases.

Case 1. $y_1 > 0$. In this case, we take a special direction in the limit (4.3) for $u \to \theta$ by $u_1 = 0$ and $u_2 \downarrow 0$. It follows that

$$\limsup_{u\to\theta} \frac{x_1 u_1 + x_2 u_2 - y_1 \frac{u_1^2-u_2^2}{u_1^2+u_2^2} - y_2 \frac{2u_1 u_2}{u_1^2+u_2^2}}{\|u\|+1}$$

$$\geq \lim_{u_1=0 \text{ and } u_2\downarrow 0} \frac{x_2 u_2 + y_1}{u_2+1}$$

$$= y_1 > 0.$$

This implies that if $y_1 > 0$, then $x \notin \widehat{D}^*f(\theta)(y)$.

Case 2. $y_1 < 0$. In this case, we take a special direction in the limit (4.3) for $u \to \theta$ by $u_2 = 0$ and $u_1 \downarrow 0$. We have

$$\limsup_{u\to\theta} \frac{x_1 u_1 + x_2 u_2 - y_1 \frac{u_1^2-u_2^2}{u_1^2+u_2^2} - y_2 \frac{2u_1 u_2}{u_1^2+u_2^2}}{\|u\|+\|f(u)\|}$$

$$\geq \lim_{u_2=0 \text{ and } u_1\downarrow 0} \frac{x_1 u_1 - y_1}{u_1+1}$$

$$= -y_1 > 0.$$

This implies that if $y_1 < 0$, then $x \notin \widehat{D}^* f(\theta)(y)$.

Case 3. $y_2 > 0$. In this case, we take a special direction in the limit (4.3) for $u \to \theta$ by $u_1 = -u_2$ and $u_1 \downarrow 0$. We have

$$\limsup_{u \to \theta} \frac{x_1 u_1 + x_2 u_2 - y_1 \frac{u_1^2 - u_2^2}{u_1^2 + u_2^2} - y_2 \frac{2 u_1 u_2}{u_1^2 + u_2^2}}{\|u\| + \|f(u)\|}$$

$$\geq \lim_{u_1 = -u_2 \text{ and } u_1 \downarrow 0} \frac{x_1 u_1 - x_2 u_1 + \frac{2 y_2}{2}}{\sqrt{2} u_1 + 1}$$

$$= y_2 > 0.$$

This implies that if $y_2 > 0$, then $x \notin \widehat{D}^* f(\theta)(y)$.

Case 4. $y_2 < 0$. In this case, we take a special direction in the limit (4.3) for $u \to \theta$ by $u_1 = u_2$ and $u_1 \downarrow 0$. We have

$$\limsup_{u \to \theta} \frac{x_1 u_1 + x_2 u_2 - y_1 \frac{u_1^2 - u_2^2}{u_1^2 + u_2^2} - y_2 \frac{2 u_1 u_2}{u_1^2 + u_2^2}}{\|u\| + \|f(u)\|}$$

$$\geq \lim_{u_1 = u_2 \text{ and } u_1 \uparrow 0} \frac{x_1 u_1 + x_2 u_2 - \frac{2 y_2}{2}}{\sqrt{2} |u_1| + 1}$$

$$= - y_2$$

$$> 0.$$

This implies that if $y_2 < 0$, then $x \notin \widehat{D}^* f(\theta)(y)$. Summarizing the above 4 cases, we obtain that

$$y \neq \theta \implies x \notin \widehat{D}^* f(\theta)(y), \text{ for any } x \in \mathbb{R}^2.$$

This implies that $\widehat{D}^* f(\theta)(y) = \emptyset$, for any $y \neq \theta$, which proves (d). □

**A trigonometric mapping 4.2.** We define $f: \mathbb{R}^2 \to \mathbb{R}^2$ by

$$f(x_1, x_2) = (\sin(x_1 + x_2), \cos(x_1 + x_2)), \text{ for any } (x_1, x_2) \in \mathbb{R}^2.$$

Then, $f$ has the following properties.

(a) $f$ is a norm constant mapping with

$$\|f(x)\| = 1, \text{ for any } x \in \mathbb{R}^2.$$

(b) $f$ is Fréchet and Mordukhovich differentiable on $\mathbb{R}^2$. For any $z = (z_1, z_2) \in \mathbb{R}^2$, we have

$$\nabla f(z) = \begin{pmatrix} \cos(z_1 + z_2) & -\sin(z_1 + z_2) \\ \cos(z_1 + z_2) & -\sin(z_1 + z_2) \end{pmatrix},$$

and

$$\widehat{D}^* f(z) = \begin{pmatrix} \cos(z_1 + z_2) & \cos(z_1 + z_2) \\ -\sin(z_1 + z_2) & -\sin(z_1 + z_2) \end{pmatrix}.$$

(c) The covering constant for $f$ is constant with

$$\hat{\alpha}(f,\bar{z},\bar{w}) = 0, \text{ for any } \bar{z} = (\bar{z}_1,\bar{z}_2) \in \mathbb{R}^2 \text{ with } \bar{w} = f(\bar{z}).$$

*Proof.* The proofs of (a) and (b) are straight forward and they are omitted here. We only prove (c).

Let $x = (x_1, x_2)$ and $y = (y_1, y_2) \in \mathbb{R}^2$. If $x = \hat{D}^* f(z)(y)$, then

$$x_1 = y_1 \cos(z_1 + z_2) - y_2 \sin(z_1 + z_2),$$

$$x_2 = y_1 \cos(z_1 + z_2) - y_2 \sin(z_1 + z_2).$$

This implies that

$$x_1^2 + x_2^2 = (y_1 \cos(z_1 + z_2) - y_2 \sin(z_1 + z_2))^2 + (y_1 \cos(z_1 + z_2) - y_2 \sin(z_1 + z_2))^2$$

$$= 2(y_1 \cos(z_1 + z_2) - y_2 \sin(z_1 + z_2))^2.$$

Then we obtain that

$$\|x\| = \sqrt{2}|y_1 \cos(z_1 + z_2) - y_2 \sin(z_1 + z_2)|, \text{ for } x = \hat{D}^* f(z)(y). \tag{4.4}$$

For any given point $\bar{z} = (\bar{z}_1, \bar{z}_2) \neq \theta$ and $\bar{w} = (\bar{w}_1, \bar{w}_2) \in \mathbb{R}^2$ with $\bar{w} = f(\bar{z})$, by (4.4), we calculate

$$\hat{\alpha}(f,\bar{z},\bar{w}) = \sup_{\eta>0} \inf\{\|x\|: x \in \hat{D}^* f(z)(y), z \in \mathbb{B}(\bar{z},\eta), f(z) \in \mathbb{B}(\bar{w},\eta), \|y\| = 1\}.$$

By the continuity of $f$ on $\mathbb{R}^2$, for any given $\eta > 0$, we have

$$\{z \in \mathbb{B}(\bar{z},\eta): f(z) \in \mathbb{B}(\bar{w},\eta)\} \neq \emptyset.$$

For any given $\eta > 0$, and $z = (z_1, z_2) \in \mathbb{B}(\bar{z}, \eta)$ and $f(z) \in \mathbb{B}(\bar{w}, \eta)$, the following system of equations has solutions with respect to $y = (y_1, y_2)$:

$$\begin{cases} y_1^2 + y_2^2 = 1, \\ y_1 \cos(z_1 + z_2) - y_2 \sin(z_1 + z_2) = 0. \end{cases} \tag{4.5}$$

By (4.4) and (4.5), we have

$$\hat{\alpha}(f,\bar{z},\bar{w}) = \sup_{\eta>0} \inf\{\|x\|: x \in \hat{D}^* f(z)(y), z \in \mathbb{B}(\bar{z},\eta), f(z) \in \mathbb{B}(\bar{w},\eta), \|y\| = 1\}$$

$$\leq \sup_{\eta>0} \inf\{\|x\|: x \in \hat{D}^* f(z)(y), z \in \mathbb{B}(\bar{z},\eta)\setminus\{\theta\}, f(z) \in \mathbb{B}(\bar{w},\eta), \|y\| = 1, y_1 \cos(z_1 + z_2) - y_2 \sin(z_1 + z_2) = 0\}$$

$$= \sup_{\eta>0} \inf\{0: x \in \hat{D}^* f(z)(y), z \in \mathbb{B}(\bar{z},\eta)\setminus\{\theta\}, f(z) \in \mathbb{B}(\bar{w},\eta), \|y\| = 1, y_1 \cos(z_1 + z_2) - y_2 \sin(z_1 + z_2) = 0\}$$

$$= 0.$$

This proves that

$$\hat{\alpha}(f,\bar{z},\bar{w}) = 0, \text{ for any } \bar{z} = (\bar{z}_1, \bar{z}_2) \in \mathbb{R}^2 \text{ with } \bar{w} = f(\bar{z}).$$

**A polynomial mapping 4.3.** We define $f: \mathbb{R}^2 \to \mathbb{R}^2$ by

$$f(x_1, x_2) = (x_1^2 - x_2^2, 2x_1 x_2), \text{ for any } (x_1, x_2) \in \mathbb{R}^2.$$

Then, $f$ has the following properties.

(a) $f$ is a norm-expansion mapping with

$$\|f(x)\| = \|x\|^2, \text{ for any } x \in \mathbb{R}^2.$$

(b) $f$ is Fréchet differentiable and Mordukhovich differentiable $\mathbb{R}^2$. For any $z = (z_1, z_2) \in \mathbb{R}^2$,

$$\nabla f(z) = \begin{pmatrix} 2z_1 & 2z_2 \\ -2z_2 & 2z_1 \end{pmatrix},$$

and

$$\widehat{D}^* f(z) = \begin{pmatrix} 2z_1 & -2z_2 \\ 2z_2 & 2z_1 \end{pmatrix}.$$

(c) The covering constant for $f$ satisfies

$$\hat{\alpha}(f, \bar{z}, f(\bar{z})) = 2\|\bar{z}\|, \text{ for any } \bar{z} \in \mathbb{R}^2.$$

In particular, we have $\hat{\alpha}(f, \theta, \theta) = 0$.

*Proof.* Proof of (a). This is a polynomial mapping. We calculate the norms of this mapping.

$$\|f(x_1, x_2)\|^2 = x_1^4 - 2x_1^2 x_2^2 + x_2^4 + 4x_1^2 x_2^2 = (x_1^2 + x_2^2)^2 = \|(x_1, x_2)\|^4.$$

This implies that $f$ is a norm-expansion mapping with

$$\|f(x)\| = \|x\|^2, \text{ for any } x \in \mathbb{R}^2.$$

Proof of (b). Since $f$ is a polynomial mapping. By Proposition 3.3, $f$ is Fréchet differentiable and Mordukhovich differentiable at every point $z = (z_1, z_2) \in \mathbb{R}^2$. By Theorem 3.1, we obtain $\nabla f(z)$ and $\widehat{D}^* f(z)$ as given in Part (b).

Proof of (c). Let $z = (z_1, z_2)$, $x = (x_1, x_2)$ and $y = (y_1, y_2) \in \mathbb{R}^2$. If $x = \widehat{D}^* f(z)(y)$, we have that

$$x_1 = 2y_1 z_1 + 2y_2 z_2,$$

$$x_2 = -2y_1 z_2 + 2y_2 z_1.$$

We calculate

$$x_1^2 + x_2^2 = (2y_1 z_1 + 2y_2 z_2)^2 + (-2y_1 z_2 + 2y_2 z_1)^2$$

$$= 4y_1^2 z_1^2 + 8y_1 y_2 z_1 z_2 + 4y_2^2 z_2^2 + 4y_1^2 z_2^2 - 8y_1 y_2 z_1 z_2 + 4y_2^2 z_1^2$$

$$= 4(y_1^2 + y_2^2)(z_1^2 + z_2^2).$$

This implies

$$\|x\| = 2\|y\|\|z\|, \text{ for } x = \widehat{D}^* f(z)(y). \tag{4.5}$$

For any given point $\bar{z} = (\bar{z}_1, \bar{z}_2)$ and $\bar{w} \in \mathbb{R}^2$ with $\bar{w} = f(\bar{z})$, by the continuity of $f$ on $\mathbb{R}^2$, we have that,

$$\{z \in \mathbb{B}(\bar{z}, \eta) : f(z) \in \mathbb{B}(\bar{w}, \eta)\} \neq \emptyset, \text{ for any } \eta > 0. \tag{4.6}$$

By the Mordukhovich differentiability of $f$ on $\mathbb{R}^2$ in Part (a) of this example, we have

$$\widehat{D}^* f(z)(y) \neq \emptyset, \text{ for any } y \in \mathbb{R}^2. \tag{4.7}$$

By (4.6) and (4.7), we obtain that, for any $\eta > 0$,

$$\{x = \widehat{D}^* f(z)(y) : z \in \mathbb{B}(\bar{z}, \eta), f(z) \in \mathbb{B}(\bar{w}, \eta), \|y\| = 1\} \neq \emptyset. \tag{4.8}$$

For any positive numbers $p$ and $q$ with $p < q$, we have

$$\{z \in \mathbb{B}(\bar{z}, p) : f(z) \in \mathbb{B}(\bar{w}, p)\} \subseteq \{z \in \mathbb{B}(\bar{z}, q) : f(z) \in \mathbb{B}(\bar{w}, q)\}. \tag{4.9}$$

By (4.6) to (4.9), for any $0 < p < q$, we get that

$$\{x = \widehat{D}^* f(z)(y) : z \in \mathbb{B}(\bar{z}, p), f(z) \in \mathbb{B}(\bar{w}, p), \|y\| = 1\}$$

$$\subseteq \{x = \widehat{D}^* f(z)(y) : z \in \mathbb{B}(\bar{z}, q), f(z) \in \mathbb{B}(\bar{w}, q), \|y\| = 1\}.$$

This implies that $\{x = \widehat{D}^* f(z)(y) : z \in \mathbb{B}(\bar{z}, \eta), f(z) \in \mathbb{B}(\bar{w}, \eta), \|y\| = 1\}$ is a decreasing (inclusion) net with respect to $\eta \downarrow 0$. Hence, the following net

$$\inf\{\|x\| : x = \widehat{D}^* f(z)(y), z \in \mathbb{B}(\bar{z}, \eta), f(z) \in \mathbb{B}(\bar{w}, \eta), \|y\| = 1\}, \tag{4.10}$$

is an increasing net of nonnegative numbers. By (4.5), we estimate (4.10).

$$\inf\{\|x\| : x = \widehat{D}^* f(z)(y), z \in \mathbb{B}(\bar{z}, \eta), f(z) \in \mathbb{B}(\bar{w}, \eta), \|y\| = 1\}$$

$$\leq \inf\{\|x\| : x = \widehat{D}^* f(\bar{z})(y), \bar{z} \in \mathbb{B}(\bar{z}, \eta), f(\bar{z}) = \bar{w} \in \mathbb{B}(\bar{w}, \eta), \|y\| = 1\}$$

$$= \inf\{2\|\bar{z}\| : x = \widehat{D}^* f(\bar{z})(y), \bar{z} \in \mathbb{B}(\bar{z}, \eta), f(\bar{z}) = \bar{w} \in \mathbb{B}(\bar{w}, \eta), \|y\| = 1\}$$

$$= 2\|\bar{z}\|.$$

This implies that (4.10) is a bounded increasing net of nonnegative numbers. We calculate

$$\hat{\alpha}(f, \bar{z}, f(\bar{z})) = \sup_{\eta > 0} \inf\{\|x\| : x = \widehat{D}^* f(z)(y), z \in \mathbb{B}(\bar{z}, \eta), f(z) \in \mathbb{B}(\bar{w}, \eta), \|y\| = 1\}$$

$$= \lim_{\eta \downarrow 0} \inf\{2\|y\| \|z\|, z \in \mathbb{B}(\bar{z}, \eta), f(z) \in \mathbb{B}(\bar{w}, \eta), \|y\| = 1\}$$

$$= \lim_{\eta \downarrow 0} \inf\{2\|z\|, z \in \mathbb{B}(\bar{z}, \eta), f(z) \in \mathbb{B}(\bar{w}, \eta), \|y\| = 1\}$$

$$= \lim_{\eta \downarrow 0} \inf\{2\|z\|, z \in \mathbb{B}(\bar{z}, \eta), f(z) \in \mathbb{B}(\bar{w}, \eta), \|y\| = 1\}$$

$$= 2\|\bar{z}\|. \qquad \square$$

**An exponential mapping 4.4.** We define $f : \mathbb{R}^2 \to \mathbb{R}^2$ by

$$f(x) = (e^{x_1+x_2},\ e^{-x_1-x_2}), \text{ for any } x = (x_1, x_2) \in \mathbb{R}^2.$$

Then, $f$ is a continuous mapping on $\mathbb{R}^2$, which has the following properties.

(a) $f$ is Fréchet differentiable and Mordukhovich differentiable on $\mathbb{R}^2$. For each $z = (z_1, z_2) \in \mathbb{R}^2$, we have

$$\nabla f(z) = \begin{pmatrix} e^{z_1+z_2} & -e^{-z_1-z_2} \\ e^{z_1+z_2} & -e^{-z_1-z_2} \end{pmatrix},$$

and

$$\widehat{D}^* f(z) = \begin{pmatrix} e^{z_1+z_2} & e^{z_1+z_2} \\ -e^{-z_1-z_2} & -e^{-z_1-z_2} \end{pmatrix}.$$

(b) The covering constant for $f$ is constant satisfying

$$\hat{\alpha}(f, \bar{z}, \bar{w}) = 0, \text{ for any } \bar{z} = (\bar{z}_1, \bar{z}_2) \in \mathbb{R}^2 \text{ with } \bar{w} = f(\bar{z}).$$

*Proof.* The proof of (a) is straight forward and it is omitted here. We prove (b). Let $\bar{z} = (\bar{z}_1, \bar{z}_2) \in \mathbb{R}^2$ with $\bar{w} = f(\bar{z})$. Let $x = (x_1, x_2)$ and $y = (y_1, y_2) \in \mathbb{R}^2$, if $x = \widehat{D}^* f(z)(y)$, by part (a), we have that

$$x_1 = y_1 e^{z_1+z_2} - y_2 e^{-z_1-z_2},$$

$$x_2 = y_1 e^{z_1+z_2} - y_2 e^{-z_1-z_2}.$$

Hence, if $x = \widehat{D}^* f(z)(y)$, then

$$x_1^2 + x_2^2 = (y_1 e^{z_1+z_2} - y_2 e^{-z_1-z_2})^2 + (y_1 e^{z_1+z_2} - y_2 e^{-z_1-z_2})^2$$

$$= 2(y_1^2 e^{2(z_1+z_2)} - 2y_1 y_2 + y_2^2 e^{-2(z_1+z_2)}).$$

This implies that

$$\|x\| = \sqrt{2}|y_1 e^{z_1+z_2} - y_2 e^{-z_1-z_2}|$$

$$= \sqrt{2}\sqrt{y_1^2 e^{2(z_1+z_2)} - 2y_1 y_2 + y_2^2 e^{-2(z_1+z_2)}}, \text{ for } x = \widehat{D}^* f(z)(y). \quad (4.11)$$

For any given point $\bar{z} = (\bar{z}_1, \bar{z}_2)$ and $\bar{w} = (\bar{w}_1, \bar{w}_2) \in \mathbb{R}^2$ with $\bar{w} = f(\bar{z})$, we calculate

$$\hat{\alpha}(f, \bar{z}, \bar{w}) = \sup_{\eta > 0} \inf\{\|x\| : x \in \widehat{D}^* f(z)(y), z \in \mathbb{B}(\bar{z}, \eta), f(z) \in \mathbb{B}(\bar{w}, \eta), \|y\| = 1\}.$$

By the continuity of $f$ on $\mathbb{R}^2$, for any given $\eta > 0$, we have

$$\{z \in \mathbb{B}(\bar{z}, \eta) : f(z) \in \mathbb{B}(\bar{w}, \eta)\} \neq \emptyset.$$

For any given $\eta > 0$, and $z = (z_1, z_2) \in \mathbb{B}(\bar{z}, \eta)$ and $f(z) \in \mathbb{B}(\bar{w}, \eta)$, one can prove that the following system of equations has solutions with respect to $y = (y_1, y_2)$:

$$\begin{cases} y_1^2 + y_2^2 = 1, \\ y_1 e^{z_1+z_2} - y_2 e^{-z_1-z_2} = 0. \end{cases}$$

By (4.11), this implies that

$$\hat{\alpha}(f,\bar{z},\bar{w}) = \sup_{\eta>0} \inf\{\|x\|: x \in \widehat{D}^*f(z)(y), z \in \mathbb{B}(\bar{z},\eta), f(z) \in \mathbb{B}(\bar{w},\eta), \|y\| = 1\}$$

$$\leq \sup_{\eta>0} \inf\{\|x\|: x \in \widehat{D}^*f(z)(y), z = (z_1,z_2) \in \mathbb{B}(\bar{z},\eta), f(z) \in \mathbb{B}(\bar{w},\eta), \|y\| = 1, y_1 e^{z_1+z_2} - y_2 e^{-z_1-z_2} = 0\}$$

$$\leq \sup_{\eta>0} \inf\{\sqrt{2}|y_1 e^{z_1+z_2} - y_2 e^{-z_1-z_2}|: x \in \widehat{D}^*f(z)(y), z \in \mathbb{B}(\bar{z},\eta), f(z) \in \mathbb{B}(\bar{w},\eta),$$
$$\|y\| = 1, y_1 e^{z_1+z_2} - y_2 e^{-z_1-z_2} = 0\}$$

$$= \sup_{\eta>0} \inf\{0: x \in \widehat{D}^*f(z)(y), z \in \mathbb{B}(\bar{z},\eta)\setminus\{\theta\}, f(z) \in \mathbb{B}(\bar{w},\eta), \|y\| = 1, y_1 e^{z_1+z_2} - y_2 e^{-z_1-z_2} = 0\}$$

$$= 0.$$

This proves that

$$\hat{\alpha}(f,\bar{z},\bar{w}) = 0, \text{ for any } \bar{z} = (\bar{z}_1,\bar{z}_2) \in \mathbb{R}^2 \text{ with } \bar{w} = f(\bar{z}). \qquad \square$$

**A logarithmic mapping 4.5.** We define $f: \mathbb{R}^2 \to \mathbb{R}^2$ by

$$f(x) = \left(\ln(1 + x_1^2 + x_2^2), \tfrac{1}{1+x_1^2+x_2^2}\right), \text{ for } x = (x_1,x_2) \in \mathbb{R}^2.$$

Then, $f$ is a continuous mapping on $\mathbb{R}^2$, which has the following properties.

(a) $f$ is Fréchet differentiable and Mordukhovich differentiable on $\mathbb{R}^2$. For any $z = (z_1,z_2) \in \mathbb{R}^2$, we have

$$\nabla f(z) = \begin{pmatrix} \frac{2z_1}{1+z_1^2+z_2^2} & -\frac{2z_1}{(1+z_1^2+z_2^2)^2} \\ \frac{2z_2}{1+z_1^2+z_2^2} & -\frac{2z_2}{(1+z_1^2+z_2^2)^2} \end{pmatrix},$$

and

$$\widehat{D}^*f(z) = \begin{pmatrix} \frac{2z_1}{1+z_1^2+z_2^2} & \frac{2z_2}{1+z_1^2+z_2^2} \\ -\frac{2z_1}{(1+z_1^2+z_2^2)^2} & -\frac{2z_2}{(1+z_1^2+z_2^2)^2} \end{pmatrix}.$$

(b) The covering constant for $f$ is constant on $\mathbb{R}^2$ satisfying

$$\hat{\alpha}(f,\bar{z},\bar{w}) = 0, \text{ for any } \bar{z} = (\bar{z}_1,\bar{z}_2) \text{ with } \bar{w} = f(\bar{z}).$$

*Proof.* By Lemma 3.1 and Theorem 3.2, the proof of (a) is straight forward and it is omitted here. We prove (b). Let $\bar{z} = (\bar{z}_1,\bar{z}_2)$ with $\bar{w} = f(\bar{z})$. Let $x = (x_1,x_2)$ and $y = (y_1,y_2) \in \mathbb{R}^2$, if $x = \widehat{D}^*f(z)(y)$, by part (a), we have that

$$x_1 = y_1 \frac{2z_1}{1+z_1^2+z_2^2} - y_2 \frac{2z_1}{(1+z_1^2+z_2^2)^2},$$

$$x_2 = y_1 \frac{2z_2}{1+z_1^2+z_2^2} - y_2 \frac{2z_2}{(1+z_1^2+z_2^2)^2}.$$

This implies that

$$x_1^2 + x_2^2 = \left(y_1 \frac{2z_1}{1+z_1^2+z_2^2} - y_2 \frac{2z_1}{(1+z_1^2+z_2^2)^2}\right)^2 + \left(y_1 \frac{2z_2}{1+z_1^2+z_2^2} - y_2 \frac{2z_2}{(1+z_1^2+z_2^2)^2}\right)^2$$

$$= y_1^2 \frac{4z_1^2}{(1+z_1^2+z_2^2)^2} - y_1 y_2 \frac{8z_1^2}{(1+z_1^2+z_2^2)^3} + y_2^2 \frac{4z_1^2}{(1+z_1^2+z_2^2)^4} + y_1^2 \frac{4z_2^2}{(1+z_1^2+z_2^2)^2} - y_1 y_2 \frac{8z_2^2}{(1+z_1^2+z_2^2)^3} + y_2^2 \frac{4z_2^2}{(1+z_1^2+z_2^2)^4}$$

$$= y_1^2 \frac{4(z_1^2+z_2^2)}{(1+z_1^2+z_2^2)^2} - y_1 y_2 \frac{8(z_1^2+z_2^2)}{(1+z_1^2+z_2^2)^3} + y_2^2 \frac{4(z_1^2+z_2^2)}{(1+z_1^2+z_2^2)^4}$$

$$= \frac{4(z_1^2+z_2^2)}{(1+z_1^2+z_2^2)^2} \left(y_1 - \frac{y_2}{1+z_1^2+z_2^2}\right)^2.$$

This implies that

$$\|x\| = \frac{2\sqrt{z_1^2+z_2^2}}{1+z_1^2+z_2^2} \left|y_1 - \frac{y_2}{1+z_1^2+z_2^2}\right|, \text{ for } x = \widehat{D}^* f(z)(y).$$

Then, we calculate

$$\hat{\alpha}(f, \bar{z}, \bar{w}) = \sup_{\eta > 0} \inf\{\|x\| : x \in \widehat{D}^* f(z)(y), z \in \mathbb{B}(\bar{z}, \eta), f(z) \in \mathbb{B}(\bar{w}, \eta), \|y\| = 1\}$$

$$= \sup_{\eta > 0} \inf\left\{\frac{2\sqrt{z_1^2+z_2^2}}{1+z_1^2+z_2^2} \left|y_1 - \frac{y_2}{1+z_1^2+z_2^2}\right| : x \in \widehat{D}^* f(z)(y), z \in \mathbb{B}(\bar{z}, \eta), f(z) \in \mathbb{B}(\bar{w}, \eta), \|y\| = 1\right\}$$

$$\leq \sup_{\eta > 0} \inf\left\{0 : x \in \widehat{D}^* f(z)(y), z \in \mathbb{B}(\bar{z}, \eta), f(z) \in \mathbb{B}(\bar{w}, \eta), \|y\| = 1, y_1 - \frac{y_2}{1+z_1^2+z_2^2} = 0\right\} = 0. \quad \square$$

**A rational mapping with radical 4.6.** We define $f: \mathbb{R}^2 \to \mathbb{R}^2$, for $x = (x_1, x_2) \in \mathbb{R}^2$, by

$$f(x) = \begin{cases} \left(\frac{x_1^2}{\sqrt{x_1^2+x_2^2}}, \frac{x_2^2}{\sqrt{x_1^2+x_2^2}}\right), & \text{for } (x_1, x_2) \neq \theta, \\ \theta, & \text{for } (x_1, x_2) = \theta. \end{cases}$$

$f$ is a continuous mapping on $\mathbb{R}^2$, which has the following properties.

(a) $f$ is Fréchet differentiable and Mordukhovich differentiable on $\mathbb{R}^2 \setminus \{\theta\}$. For any $z = (z_1, z_2) \in \mathbb{R}^2 \setminus \{\theta\}$, we have

$$\nabla f(z) = \begin{pmatrix} \frac{z_1 z_1^2 + 2z_1 z_2^2}{(z_1^2+z_2^2)\sqrt{x_1^2+x_2^2}} & \frac{-z_2^2 z_1}{(z_1^2+z_2^2)\sqrt{x_1^2+x_2^2}} \\ \frac{-z_1^2 z_2}{(z_1^2+z_2^2)\sqrt{x_1^2+x_2^2}} & \frac{z_2 z_2^2 + 2z_2 z_1^2}{(z_1^2+z_2^2)\sqrt{x_1^2+x_2^2}} \end{pmatrix},$$

and

$$\widehat{D}^* f(z) = \begin{pmatrix} \frac{z_1 z_1^2 + 2z_1 z_2^2}{(z_1^2+z_2^2)\sqrt{x_1^2+x_2^2}} & \frac{-z_1^2 z_2}{(z_1^2+z_2^2)\sqrt{x_1^2+x_2^2}} \\ \frac{-z_2^2 z_1}{(z_1^2+z_2^2)\sqrt{x_1^2+x_2^2}} & \frac{z_2 z_2^2 + 2z_2 z_1^2}{(z_1^2+z_2^2)\sqrt{x_1^2+x_2^2}} \end{pmatrix}.$$

(b) The covering constant for $f$ on $\mathbb{R}^2\setminus\{\theta\}$ satisfies that for any $\bar{z} = (\bar{z}_1, \bar{z}_2) \neq \theta$ with $\bar{w} = f(\bar{z})$,

(i) $\hat{\alpha}(f, \bar{z}, \bar{w}) \leq \frac{1}{\sqrt{2}}$.

(ii) $\hat{\alpha}(f, \bar{z}, \bar{w}) = 0$, if $\bar{z}_1 \bar{z}_2 = 0$.

(iii) $\hat{\alpha}(f, \bar{z}, \bar{w}) \leq \frac{2|\bar{z}_1 \bar{z}_2|}{\sqrt{\bar{z}_1^4 + \bar{z}_2^4}}$.

Notice that (ii) can be immediately reduced by (iii). Since the ideas of the direct proof of (ii) may be interesting for some readers, so we provide the direct proof of (i).

*Proof.* By Lemma 3.1 and Theorem 3.2, the proof of (a) is straight forward and it is omitted here.

Proof of (b). Let $z = (z_1, z_2) \in \mathbb{R}^2\setminus\{\theta\}$, and $x = (x_1, x_2)$, $y = (y_1, y_2) \in \mathbb{R}^2$. Suppose that $x = \widehat{D}^* f(z)(y)$. By Part (a), we have that

$$x_1 = y_1 \frac{z_1 z_1^2 + 2 z_1 z_2^2}{(z_1^2 + z_2^2)\sqrt{x_1^2 + x_2^2}} + y_2 \frac{-z_2^2 z_1}{(z_1^2 + z_2^2)\sqrt{x_1^2 + x_2^2}},$$

$$x_2 = y_1 \frac{-z_1^2 z_2}{(z_1^2 + z_2^2)\sqrt{x_1^2 + x_2^2}} + y_2 \frac{z_2 z_2^2 + 2 z_2 z_1^2}{(z_1^2 + z_2^2)\sqrt{x_1^2 + x_2^2}}.$$

This implies that

$$x_1^2 + x_2^2 = \left( y_1 \frac{z_1 z_1^2 + 2 z_1 z_2^2}{(z_1^2+z_2^2)\sqrt{x_1^2+x_2^2}} + y_2 \frac{-z_2^2 z_1}{(z_1^2+z_2^2)\sqrt{x_1^2+x_2^2}} \right)^2 + \left( y_1 \frac{-z_1^2 z_2}{(z_1^2+z_2^2)\sqrt{x_1^2+x_2^2}} + y_2 \frac{z_2 z_2^2 + 2 z_2 z_1^2}{(z_1^2+z_2^2)\sqrt{x_1^2+x_2^2}} \right)^2$$

$$= y_1^2 \frac{z_1^2 z_1^4 + 4 z_1^4 z_2^2 + 4 z_1^2 z_2^4}{(z_1^2+z_2^2)^3} - 2 y_1 y_2 \frac{z_1^4 z_2^2 + 2 z_1^2 z_2^4}{(z_1^2+z_2^2)^3} + y_2^2 \frac{z_1^2 z_2^4}{(z_1^2+z_2^2)^3}$$

$$+ y_1^2 \frac{z_1^4 z_2^2}{(z_1^2+z_2^2)^3} - 2 y_1 y_2 \frac{z_1^2 z_2^4 + 2 z_1^4 z_2^2}{(z_1^2+z_2^2)^3} + y_2^2 \frac{z_2^2 z_2^4 + 4 z_1^2 z_2^4 + 4 z_1^4 z_2^2}{(z_1^2+z_2^2)^3}$$

$$= y_1^2 \frac{z_1^2 z_1^4 + 5 z_1^4 z_2^2 + 4 z_1^2 z_2^4}{(z_1^2+z_2^2)^3} - 2 y_1 y_2 \frac{3 z_1^4 z_2^2 + 3 z_1^2 z_2^4}{(z_1^2+z_2^2)^3} + y_2^2 \frac{z_2^2 z_2^4 + 5 z_1^2 z_2^4 + 4 z_1^4 z_2^2}{(z_1^2+z_2^2)^3}$$

$$= \frac{1}{(z_1^2+z_2^2)^3} \left( y_1^2 z_1^2 (z_1^4 + 5 z_1^2 z_2^2 + 4 z_2^4) - 6 y_1 y_2 z_1^2 z_2^2 (z_1^2 + z_2^2) + y_2^2 z_2^2 (z_2^4 + 5 z_1^2 z_2^2 + 4 z_1^4) \right)$$

$$= \frac{1}{(z_1^2+z_2^2)^3} \left( y_1^2 z_1^2 (z_1^2 + 4 z_2^2)(z_1^2 + z_2^2) - 6 y_1 y_2 z_1^2 z_2^2 (z_1^2 + z_2^2) + y_2^2 z_2^2 (z_2^2 + 4 z_1^2)(z_1^2 + z_2^2) \right)$$

$$= \frac{1}{(z_1^2+z_2^2)^2} \left( y_1^2 z_1^2 (z_1^2 + 4 z_2^2) - 6 y_1 y_2 z_1^2 z_2^2 + y_2^2 z_2^2 (z_2^2 + 4 z_1^2) \right)$$

$$= \frac{1}{(z_1^2+z_2^2)^2} \left( y_1^2 z_1^4 + 4 y_1^2 z_1^2 z_2^2 - 6 y_1 y_2 z_1^2 z_2^2 + y_2^2 z_2^4 + 4 y_2^2 z_1^2 z_2^2 \right)$$

$$= \frac{1}{(z_1^2+z_2^2)^2}(y_1^2 z_1^4 + 2y_1 y_2 z_1^2 z_2^2 + y_2^2 z_2^4 + 4y_1^2 z_1^2 z_2^2 - 8y_1 y_2 z_1^2 z_2^2 + 4y_2^2 z_1^2 z_2^2)$$

$$= \frac{1}{(z_1^2+z_2^2)^2}((y_1 z_1^2 + y_2 z_2^2)^2 + 4z_1^2 z_2^2 (y_1 - y_2)^2). \tag{4.12}$$

By (4.12), for $z = (z_1, z_2) \in \mathbb{R}^2 \setminus \{\theta\}$ and $x = (x_1, x_2), y = (y_1, y_2) \in \mathbb{R}^2$, we obtain

$$\|x\| = \frac{\sqrt{(y_1 z_1^2 + y_2 z_2^2)^2 + 4z_1^2 z_2^2 (y_1 - y_2)^2}}{z_1^2 + z_2^2}, \text{ if } x = \widehat{D}^* f(z)(y).$$

In particular, if we let $y_1 = y_2 = \pm \frac{1}{\sqrt{2}}$ in (4.12), we obtain

$$\frac{1}{(z_1^2+z_2^2)^2}\left(\left(\frac{1}{\sqrt{2}} z_1^2 + \frac{1}{\sqrt{2}} z_2^2\right)^2 + 4z_1^2 z_2^2 \left(\frac{1}{\sqrt{2}} - \frac{1}{\sqrt{2}}\right)^2\right) = \frac{1}{2}. \tag{4.13}$$

For $\bar{z} = (\bar{z}_1, \bar{z}_2) \neq \theta$ with $\bar{w} = f(\bar{z})$, by (4.12) and (4.13), we have that

$$\hat{\alpha}(f, \bar{z}, \bar{w}) = \sup_{\eta > 0} \inf \{\|x\|: x \in \widehat{D}^* f(z)(y), z \in \mathbb{B}(\bar{z}, \eta) \setminus \{\theta\}, f(z) \in \mathbb{B}(\bar{w}, \eta), \|y\| = 1\}$$

$$= \sup_{\eta > 0} \inf \left\{\frac{\sqrt{(y_1 z_1^2 + y_2 z_2^2)^2 + 4z_1^2 z_2^2 (y_1 - y_2)^2}}{z_1^2 + z_2^2}: z \in \mathbb{B}(\bar{z}, \eta) \setminus \{\theta\}, f(z) \in \mathbb{B}(\bar{w}, \eta), \|y\| = 1\right\}$$

$$\leq \sup_{\eta > 0} \inf \left\{\frac{\sqrt{(y_1 z_1^2 + y_2 z_2^2)^2 + 4z_1^2 z_2^2 (y_1 - y_2)^2}}{z_1^2 + z_2^2}: z \in \mathbb{B}(\bar{z}, \eta) \setminus \{\theta\}, f(z) \in \mathbb{B}(\bar{w}, \eta), \|y\| = 1, y_1 = y_2 = \pm \frac{1}{\sqrt{2}}\right\}$$

$$= \sup_{\eta > 0} \inf \left\{\frac{1}{\sqrt{2}}: z \in \mathbb{B}(\bar{z}, \eta) \setminus \{\theta\}, f(z) \in \mathbb{B}(\bar{w}, \eta), \|y\| = 1, y_1 = y_2 = \pm \frac{1}{\sqrt{2}}\right\}$$

$$= \frac{1}{\sqrt{2}}.$$

This proves (i) in Part (b). Next, we suppose $\bar{z}_1 \bar{z}_2 = 0$ with $\bar{z}_i = 0$, for some $k = 1, 2$. By assumption that $\bar{z} = (\bar{z}_1, \bar{z}_2) \neq \theta$, we have $\bar{z}_{3-k} \neq 0$. By the continuity of $f$ around $\bar{z} = (\bar{z}_1, \bar{z}_2) \neq \theta$, for any $\eta > 0$, we can have $z = (z_1, z_2) \in \mathbb{B}(\bar{z}, \eta) \setminus \{\theta\}$ satisfying $f(z) \in \mathbb{B}(\bar{w}, \eta)$ and $z_k = 0, z_{3-k} \neq 0$.

$$\hat{\alpha}(f, \bar{z}, \bar{w}) = \sup_{\eta > 0} \inf \{\|x\|: x \in \widehat{D}^* f(z)(y), z \in \mathbb{B}(\bar{z}, \eta) \setminus \{\theta\}, f(z) \in \mathbb{B}(\bar{w}, \eta), \|y\| = 1\}$$

$$= \sup_{\eta > 0} \inf \left\{\frac{\sqrt{(y_1 z_1^2 + y_2 z_2^2)^2 + 4z_1^2 z_2^2 (y_1 - y_2)^2}}{z_1^2 + z_2^2}: z \in \mathbb{B}(\bar{z}, \eta) \setminus \{\theta\}, f(z) \in \mathbb{B}(\bar{w}, \eta), \|y\| = 1\right\}$$

$$\leq \sup_{\eta > 0} \inf \left\{\frac{\sqrt{(y_1 z_1^2 + y_2 z_2^2)^2 + 4z_1^2 z_2^2 (y_1 - y_2)^2}}{z_1^2 + z_2^2}: z \in \mathbb{B}(\bar{z}, \eta) \setminus \{\theta\}, f(z) \in \mathbb{B}(\bar{w}, \eta), z_k = 0, z_{3-k} \neq 0, \|y\| = 1\right\}$$

$$\leq \sup_{\eta > 0} \inf \left\{\frac{|y_k 0 + 0 z_{3-k}^2|}{z_1^2 + z_2^2}: z \in \mathbb{B}(\bar{z}, \eta) \setminus \{\theta\}, f(z) \in \mathbb{B}(\bar{w}, \eta), z_k = 0, z_{3-k} \neq 0, \|y\| = 1, y_k = 1, y_{3-k} = 0\right\}$$

$$\leq \sup_{\eta>0} \inf\{0 : z \in \mathbb{B}(\bar{z},\eta)\setminus\{\theta\}, f(z) \in \mathbb{B}(\bar{w},\eta), z_k = 0, z_{3-k} \neq 0. \|y\| = 1, y_k = 1, y_{3-k} = 0\}$$

$$= 0.$$

This proves (ii) in Part (b). Next, we prove (iii) in (b). By (4.12), we estimate

$$\hat{\alpha}(f,\bar{z},\bar{w}) = \sup_{\eta>0} \inf\{\|x\| : x \in \widehat{D}^*f(z)(y), z \in \mathbb{B}(\bar{z},\eta)\setminus\{\theta\}, f(z) \in \mathbb{B}(\bar{w},\eta), \|y\| = 1\}$$

$$= \liminf_{\eta\downarrow 0} \left\{ \frac{\sqrt{(y_1 z_1^2 + y_2 z_2^2)^2 + 4 z_1^2 z_2^2 (y_1 - y_2)^2}}{z_1^2 + z_2^2} : z \in \mathbb{B}(\bar{z},\eta)\setminus\{\theta\}, f(z) \in \mathbb{B}(\bar{w},\eta), \|y\| = 1 \right\}$$

$$\leq \liminf_{\eta\downarrow 0} \left\{ \frac{\sqrt{(y_1 z_1^2 + y_2 z_2^2)^2 + 4 z_1^2 z_2^2 (y_1 - y_2)^2}}{z_1^2 + z_2^2} : z \in \mathbb{B}(\bar{z},\eta)\setminus\{\theta\}, f(z) \in \mathbb{B}(\bar{w},\eta), y_1 = -\frac{z_2^2}{\sqrt{z_2^4+z_2^4}}, y_2 = \frac{z_1^2}{\sqrt{z_2^4+z_2^4}} \right\}$$

$$\leq \liminf_{\eta\downarrow 0} \left\{ \frac{\sqrt{4 z_1^2 z_2^2}}{\sqrt{z_2^4+z_2^4}} : z \in \mathbb{B}(\bar{z},\eta)\setminus\{\theta\}, f(z) \in \mathbb{B}(\bar{w},\eta), y_1 = -\frac{z_2^2}{\sqrt{z_2^4+z_2^4}}, y_2 = \frac{z_1^2}{\sqrt{z_2^4+z_2^4}} \right\}$$

$$\leq \liminf_{\eta\downarrow 0} \left\{ \frac{\sqrt{4 \bar{z}_1^2 \bar{z}_2^2}}{\sqrt{\bar{z}_2^4+\bar{z}_2^4}} : \bar{z} \in \mathbb{B}(\bar{z},\eta)\setminus\{\theta\}, f(\bar{z}) \in \mathbb{B}(\bar{w},\eta), y_1 = -\frac{\bar{z}_2^2}{\sqrt{\bar{z}_2^4+\bar{z}_2^4}}, y_2 = \frac{\bar{z}_1^2}{\sqrt{\bar{z}_2^4+\bar{z}_2^4}} \right\}$$

$$= \frac{2|\bar{z}_1 \bar{z}_2|}{\sqrt{\bar{z}_2^4+\bar{z}_2^4}}.$$

This proves (iii) in (b). □

**4.7. Remarks on the rational mapping with radical 4.6.** The results (i), (ii) and (iii) of the rational mapping 4.6 only provide some special estimations for the covering constant for the considered mapping $f$ defined by (4.1). In general, to find $\hat{\alpha}(f,\bar{z},\bar{w})$, we have to find

$$\inf\left\{ \frac{\sqrt{(y_1 \bar{z}_1^2 + y_2 \bar{z}_2^2)^2 + 4 \bar{z}_1^2 \bar{z}_2^2 (y_1 - y_2)^2}}{\bar{z}_1^2 + \bar{z}_2^2} : \|y\| = 1 \right\}. \qquad (4.14)$$

It is clear that we have to solve the following minimization problem with respect to $y_1, y_2$,

$$\text{Minimizing:} \quad \frac{(y_1 \bar{z}_1^2 + y_2 \bar{z}_2^2)^2 + 4 \bar{z}_1^2 \bar{z}_2^2 (y_1 - y_2)^2}{(\bar{z}_1^2 + \bar{z}_2^2)^2},$$

$$\text{Subject to:} \quad y_1^2 + y_2^2 - 1 = 0.$$

It is equivalent to minimizing the following minimization problem

$$\text{Minimizing:} \quad (y_1 \bar{z}_1^2 + y_2 \bar{z}_2^2)^2 + 4 \bar{z}_1^2 \bar{z}_2^2 (y_1 - y_2)^2$$

$$\text{Subject to:} \quad y_1^2 + y_2^2 - 1 = 0.$$

We use the method of Lagrange multipliers to solve the above minimization problem with variables $y_1$ and $y_2$. Let

$$L(y_1, y_2) = (y_1 \bar{z}_1^2 + y_2 \bar{z}_2^2)^2 + 4\bar{z}_1^2 \bar{z}_2^2 (y_1 - y_2)^2 - \lambda(y_1^2 + y_2^2 - 1).$$

We have that

$$L_{y_1}(y_1, y_2) = 2\bar{z}_1^2 (y_1 \bar{z}_1^2 + y_2 \bar{z}_2^2) + 8\bar{z}_1^2 \bar{z}_2^2 (y_1 - y_2) - 2\lambda y_1 = 0$$

$$L_{y_2}(y_1, y_2) = 2\bar{z}_2^2 (y_1 \bar{z}_1^2 + y_2 \bar{z}_2^2) - 8\bar{z}_1^2 \bar{z}_2^2 (y_1 - y_2) - 2\lambda y_2 = 0$$

$$L_\lambda(y_1, y_2) = -y_1^2 - y_2^2 + 1 = 0.$$

We can find the solution of the above system of equations as below.

$$\lambda = \frac{\bar{z}_1^4 + \bar{z}_2^4 \pm \sqrt{(\bar{z}_1^4 + \bar{z}_2^4)^2 + 32\bar{z}_1^4 \bar{z}_2^4}}{2},$$

$$\frac{y_2}{y_1} = \frac{\bar{z}_1^4 - \bar{z}_2^4 \mp \sqrt{(\bar{z}_1^4 + \bar{z}_2^4)^2 + 32\bar{z}_1^4 \bar{z}_2^4}}{6\bar{z}_1^2 \bar{z}_2^2}.$$

$$y_1 = \frac{\pm 3\sqrt{2}\bar{z}_1^2 \bar{z}_2^2}{\sqrt{(\bar{z}_1^4 - \bar{z}_2^4)^2 \pm (\bar{z}_1^4 - \bar{z}_2^4)\sqrt{(\bar{z}_1^4 - \bar{z}_2^4)^2 + 36\bar{z}_1^4 \bar{z}_2^4} + 36\bar{z}_1^4 \bar{z}_2^4}}$$

$$y_2 = \frac{\mp \sqrt{(\bar{z}_1^4 - \bar{z}_2^4)^2 \pm (\bar{z}_1^4 - \bar{z}_2^4)\sqrt{(\bar{z}_1^4 - \bar{z}_2^4)^2 + 36\bar{z}_1^4 \bar{z}_2^4} + 18\bar{z}_1^4 \bar{z}_2^4}}{\sqrt{(\bar{z}_1^4 - \bar{z}_2^4)^2 \pm (\bar{z}_1^4 - \bar{z}_2^4)\sqrt{(\bar{z}_1^4 - \bar{z}_2^4)^2 + 36\bar{z}_1^4 \bar{z}_2^4} + 36\bar{z}_1^4 \bar{z}_2^4}}.$$

Here,

$$\left(\frac{y_2}{y_1}\right)^2 = \left(\frac{\bar{z}_1^4 - \bar{z}_2^4 \pm \sqrt{(\bar{z}_1^4 - \bar{z}_2^4)^2 + 36\bar{z}_1^4 \bar{z}_2^4}}{6\bar{z}_1^2 \bar{z}_2^2}\right)^2 = \frac{(\bar{z}_1^4 - \bar{z}_2^4)^2 \pm (\bar{z}_1^4 - \bar{z}_2^4)\sqrt{(\bar{z}_1^4 - \bar{z}_2^4)^2 + 36\bar{z}_1^4 \bar{z}_2^4} + 18\bar{z}_1^4 \bar{z}_2^4}{18\bar{z}_1^4 \bar{z}_2^4}.$$

When these solutions are substituted into (4.14), the solution of minimizing the value $\|x\|$, for $x = \widehat{D}^* f(z)(y)$ with $\|y\| = 1$ satisfies that

$$\|x\|^2 = \frac{(y_1 \bar{z}_1^2 + y_2 \bar{z}_2^2)^2 + 4\bar{z}_1^2 \bar{z}_2^2 (y_1 - y_2)^2}{(\bar{z}_1^2 + \bar{z}_2^2)^2}$$

$$= \frac{\left(\frac{3\sqrt{2}\bar{z}_1^2 \bar{z}_2^2}{\sqrt{(\bar{z}_1^4 - \bar{z}_2^4)^2 \pm (\bar{z}_1^4 - \bar{z}_2^4)\sqrt{(\bar{z}_1^4 - \bar{z}_2^4)^2 + 36\bar{z}_1^4 \bar{z}_2^4} + 36\bar{z}_1^4 \bar{z}_2^4}} \bar{z}_1^2 - \frac{\sqrt{(\bar{z}_1^4 - \bar{z}_2^4)^2 \pm (\bar{z}_1^4 - \bar{z}_2^4)\sqrt{(\bar{z}_1^4 - \bar{z}_2^4)^2 + 36\bar{z}_1^4 \bar{z}_2^4} + 18\bar{z}_1^4 \bar{z}_2^4}}{\sqrt{(\bar{z}_1^4 - \bar{z}_2^4)^2 \pm (\bar{z}_1^4 - \bar{z}_2^4)\sqrt{(\bar{z}_1^4 - \bar{z}_2^4)^2 + 36\bar{z}_1^4 \bar{z}_2^4} + 36\bar{z}_1^4 \bar{z}_2^4}} \bar{z}_2^2\right)^2}{(\bar{z}_1^2 + \bar{z}_2^2)^2}$$

$$+\frac{4\bar{z}_1^2\bar{z}_2^2\left(\frac{3\sqrt{2}\bar{z}_1^2\bar{z}_2^2}{\sqrt{(\bar{z}_1^4-\bar{z}_2^4)^2\pm(\bar{z}_1^4-\bar{z}_2^4)\sqrt{(\bar{z}_1^4-\bar{z}_2^4)^2+36\bar{z}_1^4\bar{z}_2^4}+36\bar{z}_1^4\bar{z}_2^4}}+\frac{\sqrt{(\bar{z}_1^4-\bar{z}_2^4)^2\pm(\bar{z}_1^4-\bar{z}_2^4)\sqrt{(\bar{z}_1^4-\bar{z}_2^4)^2+36\bar{z}_1^4\bar{z}_2^4}+18\bar{z}_1^4\bar{z}_2^4}}{\sqrt{(\bar{z}_1^4-\bar{z}_2^4)^2\pm(\bar{z}_1^4-\bar{z}_2^4)\sqrt{(\bar{z}_1^4-\bar{z}_2^4)^2+36\bar{z}_1^4\bar{z}_2^4}+36\bar{z}_1^4\bar{z}_2^4}}\right)^2}{(\bar{z}_1^2+\bar{z}_2^2)^2}$$

$$=\frac{\frac{18\bar{z}_1^8\bar{z}_2^4}{(\bar{z}_1^4-\bar{z}_2^4)^2\pm(\bar{z}_1^4-\bar{z}_2^4)\sqrt{(\bar{z}_1^4-\bar{z}_2^4)^2+36\bar{z}_1^4\bar{z}_2^4}+36\bar{z}_1^4\bar{z}_2^4}+\frac{\bar{z}_2^4\left((\bar{z}_1^4-\bar{z}_2^4)^2\pm(\bar{z}_1^4-\bar{z}_2^4)\sqrt{(\bar{z}_1^4-\bar{z}_2^4)^2+36\bar{z}_1^4\bar{z}_2^4}+18\bar{z}_1^4\bar{z}_2^4\right)}{(\bar{z}_1^4-\bar{z}_2^4)^2\pm(\bar{z}_1^4-\bar{z}_2^4)\sqrt{(\bar{z}_1^4-\bar{z}_2^4)^2+36\bar{z}_1^4\bar{z}_2^4}+36\bar{z}_1^4\bar{z}_2^4}}{(\bar{z}_1^2+\bar{z}_2^2)^2}$$

$$-\frac{\frac{6\sqrt{2}\bar{z}_1^4\bar{z}_2^2\sqrt{(\bar{z}_1^4-\bar{z}_2^4)^2\pm(\bar{z}_1^4-\bar{z}_2^4)\sqrt{(\bar{z}_1^4-\bar{z}_2^4)^2+36\bar{z}_1^4\bar{z}_2^4}+18\bar{z}_1^4\bar{z}_2^4}}{(\bar{z}_1^4-\bar{z}_2^4)^2\pm(\bar{z}_1^4-\bar{z}_2^4)\sqrt{(\bar{z}_1^4-\bar{z}_2^4)^2+36\bar{z}_1^4\bar{z}_2^4}+36\bar{z}_1^4\bar{z}_2^4}}{(\bar{z}_1^2+\bar{z}_2^2)^2}+\frac{4\bar{z}_1^2\bar{z}_2^2+\frac{24\sqrt{2}\bar{z}_1^4\bar{z}_2^2\sqrt{(\bar{z}_1^4-\bar{z}_2^4)^2\pm(\bar{z}_1^4-\bar{z}_2^4)\sqrt{(\bar{z}_1^4-\bar{z}_2^4)^2+36\bar{z}_1^4\bar{z}_2^4}+18\bar{z}_1^4\bar{z}_2^4}}{(\bar{z}_1^4-\bar{z}_2^4)^2\pm(\bar{z}_1^4-\bar{z}_2^4)\sqrt{(\bar{z}_1^4-\bar{z}_2^4)^2+36\bar{z}_1^4\bar{z}_2^4}+36\bar{z}_1^4\bar{z}_2^4}}{(\bar{z}_1^2+\bar{z}_2^2)^2}$$

$$=\frac{\frac{18\bar{z}_1^8\bar{z}_2^4}{(\bar{z}_1^4-\bar{z}_2^4)^2\pm(\bar{z}_1^4-\bar{z}_2^4)\sqrt{(\bar{z}_1^4-\bar{z}_2^4)^2+36\bar{z}_1^4\bar{z}_2^4}+36\bar{z}_1^4\bar{z}_2^4}+\frac{\bar{z}_2^4\left((\bar{z}_1^4-\bar{z}_2^4)^2\pm(\bar{z}_1^4-\bar{z}_2^4)\sqrt{(\bar{z}_1^4-\bar{z}_2^4)^2+36\bar{z}_1^4\bar{z}_2^4}+18\bar{z}_1^4\bar{z}_2^4\right)}{(\bar{z}_1^4-\bar{z}_2^4)^2\pm(\bar{z}_1^4-\bar{z}_2^4)\sqrt{(\bar{z}_1^4-\bar{z}_2^4)^2+36\bar{z}_1^4\bar{z}_2^4}+36\bar{z}_1^4\bar{z}_2^4}}{(\bar{z}_1^2+\bar{z}_2^2)^2}$$

$$+\frac{4\bar{z}_1^2\bar{z}_2^2+\frac{18\sqrt{2}\bar{z}_1^4\bar{z}_2^4\sqrt{(\bar{z}_1^4-\bar{z}_2^4)^2\pm(\bar{z}_1^4-\bar{z}_2^4)\sqrt{(\bar{z}_1^4-\bar{z}_2^4)^2+36\bar{z}_1^4\bar{z}_2^4}+18\bar{z}_1^4\bar{z}_2^4}}{(\bar{z}_1^4-\bar{z}_2^4)^2\pm(\bar{z}_1^4-\bar{z}_2^4)\sqrt{(\bar{z}_1^4-\bar{z}_2^4)^2+36\bar{z}_1^4\bar{z}_2^4}+36\bar{z}_1^4\bar{z}_2^4}}{(\bar{z}_1^2+\bar{z}_2^2)^2}.$$

We see that the solution of (4.14) is too complicated. It will be very nice to solve (4.14) with more efficient techniques to obtain have a simpler solution for (4.14).

## 5. Some Applications

Mordukhovich derivatives (coderivatives) of mappings (for both of set-valued and single-valued) have been widely applied to variational analysis. One of the important applications of the Mordukhovich derivatives is to define and to calculate the covering constants for both set-valued and single-valued mappings in Banach spaces. By the concept of covering constants, in [1], the authors Arutyunov Mordukhovich and Zhukovskiy proved an important theorem, which is called Arutyunov Mordukhovich Zhukovskiy Parameterized Coincidence Point Theorem (It is simply named as AMZ Theorem). The AMZ Theorem provides a very powerful tool in set-valued and variational analysis (see [1, 11−15]). The underlying spaces of the AMZ Theorem are Asplund spaces, which includes Euclidean spaces as special cases (See [1, 23] for more details). In this section, we use the results obtained in the previous section about the covering constants of mappings in $\mathbb{R}^2$ and by applying the AMZ Theorem to solve some parameterized equations in $\mathbb{R}^2$. We recall the AMZ Theorem.

**(AMZ Theorem)** *Let the Banach spaces X and Y in be Asplund and let P be a topological space. Let $F: X \rightrightarrows Y$ and $G(\cdot, \cdot): X \times P \rightrightarrows Y$ be set-valued mappings. Let $\bar{x} \in X$ and $\bar{y} \in Y$ with $\bar{y} \in F(\bar{x})$. Suppose that the following conditions are satisfied:*

(A1) *The multifunction $F: X \rightrightarrows Y$ is closed around $(\bar{x}, \bar{y})$.*

(A2) *There are neighborhoods $U \subset X$ of $\bar{x}$, $V \subset Y$ of $\bar{y}$, and O of $\bar{p} \in P$ as well as a number $\beta \geq 0$ such that the multifunction $G(\cdot, p): X \rightrightarrows Y$ is Lipschitz-like on U relative to V for each $p \in O$ with the uniform modulus $\beta$, while the multifunction $p \to G(\bar{x}, p)$ is lower/inner semicontinuous at $\bar{p}$.*

(A3) *The Lipschitzian modulus $\beta$ of $G(\cdot, p)$ is chosen as $\beta < \hat{\alpha}(F, \bar{x}, \bar{y})$, where $\hat{\alpha}(F, \bar{x}, \bar{y})$ is the covering constant of F around $(\bar{x}, \bar{y})$ taken from (2.5).*

*Then for each $\alpha > 0$ with $\beta < \alpha < \hat{\alpha}(F, \bar{x}, \bar{y})$, there exist a neighborhood $W \subset P$ of $\bar{p}$ and a single-valued mapping $\sigma\colon W \to X$ such that whenever $p \in W$ we have*

$$F(\sigma(p)) \cap G(\sigma(p), p) \neq \emptyset \quad \text{and} \quad \|\sigma(p) - \bar{x}\|_X \leq \frac{\mathrm{dist}(\bar{y}, G(\bar{x}, p))}{\alpha - \beta}.$$

In particular, if both F and G are single-valued mappings, then we get the following corollary of the AMZ Theorem.

**Corollary 2.2 in [12].** *Let the Banach spaces X and Y be Asplund and let P be a topological space. Let $F\colon X \to Y$ and $G(\cdot, \cdot)\colon X \times P \to Y$ be single-valued mappings. Let $\bar{x} \in X$ and $\bar{y} \in Y$ with $\bar{y} = F(\bar{x})$. Suppose that the following conditions are satisfied:*

(A1) *The mapping $F\colon X \to Y$ is continuous around $(\bar{x}, \bar{y})$.*

(A2) *There are neighborhoods $U \subset X$ of $\bar{x}$, $V \subset Y$ of $\bar{y}$, and O of $\bar{p} \in P$ as well as a number $\beta \geq 0$ such that the mapping $G(\cdot, p)\colon X \to Y$ satisfies the Lipschitz condition on U relative to V for each $p \in O$ with the uniform modulus $\beta$, while the mapping $p \to G(\bar{x}, p)$ is lower semicontinuous at $\bar{p}$.*

(A3) *The Lipschitzian modulus $\beta$ of $G(\cdot, p)$ is chosen as $\beta < \hat{\alpha}(F, \bar{x}, \bar{y})$, where $\hat{\alpha}(F, \bar{x}, \bar{y})$ is the covering constant of F around $(\bar{x}, \bar{y})$ taken from (2.5).*

*Then for each $\alpha > 0$ with $\beta < \alpha < \hat{\alpha}(F, \bar{x}, \bar{y})$, there exist a neighborhood $W \subset P$ of $\bar{p}$ and a single-valued mapping $\sigma\colon W \to X$ such that whenever $p \in W$ we have*

$$F(\sigma(p)) = G(\sigma(p), p) \quad \text{and} \quad \|\sigma(p) - \bar{x}\|_X \leq \frac{\|G(\bar{x}, p) - \bar{y}\|_Y}{\alpha - \beta}.$$

**Theorem 5.1.** *Let C be a topological space. Let $h(\cdot, \cdot) = (h_1(\cdot, \cdot), h_2(\cdot, \cdot))\colon \mathbb{R}^2 \times C \to \mathbb{R}^2$ be a single-valued mapping. Let $\omega = (\omega_1, \omega_2)\colon C \to \mathbb{R}^2$ be a single-valued lower semicontinuous mapping. Let $f\colon \mathbb{R}^2 \to \mathbb{R}^2$ be a single-valued mapping defined by*

$$f(x_1, x_2) = (x_1^2 - x_2^2, 2x_1 x_2), \text{ for any } (x_1, x_2) \in \mathbb{R}^2.$$

*Let $\bar{x}$ and $\bar{y} \in \mathbb{R}^2 \{\theta\}$ with $\bar{y} = f(\bar{x})$. Let $\bar{s} \in C$. Let $\mathbb{B}(\bar{x}, \|\bar{x}\|)$ be the closed ball in $\mathbb{R}^2$ with radius $\|\bar{x}\|$ and centered at $\bar{x}$. Suppose that the following conditions are satisfied:*

*There is a number $\beta \in [0, \|\bar{x}\|)$ such that the mapping $h(\cdot, s)\colon \mathbb{R}^2 \to \mathbb{R}^2$ satisfies the Lipschitz condition on $\mathbb{B}\left(\bar{x}, \frac{\|\bar{x}\|}{2}\right)$ for each $s \in C$ with the uniform modulus $\beta$, while the mapping $s \to h(\bar{x}, s)$ is lower semicontinuous at $\bar{s}$.*

*Then for each $\alpha > 0$ with $\beta < \alpha < \|\bar{x}\|$, there exist a neighborhood $W \subset C$ of $\bar{s}$ and a single-valued mapping $\sigma = (\sigma_1, \sigma_2)\colon W \to \mathbb{R}^2$ such that whenever $s \in W$ we have*

$$\begin{cases} \sigma_1(s)^2 - \sigma_2(s)^2 = h_1\big((\sigma_1(s), \sigma_2(s)), s\big) + \omega_1(s) \\ 2\sigma_1(s)\sigma_2(s) = h_2\big((\sigma_1(s), \sigma_2(s)), s\big) + \omega_2(s) \end{cases},$$

*and* $$\|(\sigma_1(s), \sigma_2(s)) - (\bar{x}_1, \bar{x}_2)\| \leq \frac{\|h((\bar{x}_1,\bar{x}_2),s) + (\omega_1(s),\omega_2(s)) - (\bar{x}_1^2 - \bar{x}_2^2,\ 2\bar{x}_1\bar{x}_2)\|}{\alpha - \beta}. \qquad (5.1)$$

*In particular, σ satisfies that, whenever $s \in W$,*

$$\left(h_1\big((\sigma_1(s), \sigma_2(s)), s\big) + \omega_1(s)\right)^2 + \left(h_2\big((\sigma_1(s), \sigma_2(s)), s\big) + \omega_2(s)\right)^2 = (\sigma_1(s)^2 + \sigma_2(s)^2)^2. \qquad (5.2)$$

*Proof.* For the given function $f$ in this theorem, it is studied in Example 4.3, from which we have that

$$\hat{a}(f, x, f(x)) = 2\|x\|, \text{ for any } x \in \mathbb{R}^2. \qquad (5.3)$$

Notice that, for any $x \in \mathbb{B}\left(\bar{x}, \frac{\|\bar{x}\|}{2}\right)$, we have that $\|x\| \geq \frac{\|\bar{x}\|}{2}$. By (5.3), we have that

$$\hat{a}(f, x, f(x)) = 2\|x\| \geq \|\bar{x}\|, \text{ for any } x \in \mathbb{B}\left(\bar{x}, \frac{\|\bar{x}\|}{2}\right). \qquad (5.4)$$

In AMZ Theorem, we particularly take

$$F(x) = f(x) \text{ and } G(x, s) = h(x, s) + \omega(s), \text{ for all } x \in \mathbb{R}^2, s \in C. \qquad (5.5)$$

As single-valued mappings, we can verify that the mappings, which are defined in (5.4), $F: \mathbb{R}^2 \to \mathbb{R}^2$ and $G(\cdot, \cdot): \mathbb{R}^2 \times C \to \mathbb{R}^2$ satisfy all conditions in the AMZ Theorem, from which the conclusions (5.1) are immediately proved. (5.2) is reduced by (5.1) immediately. □

## 6. Conclusion and Remarks

In this paper, we prove the guidelines for calculating the Fréchet derivatives of single-valued mappings in Euclidean spaces, which is represented by the linear approximation. With great help of the relationship between Fréchet derivatives and Mordukhovich derivatives (It is given by Theorem 1.38 in [18]), we derive the algorithm for calculating the Mordukhovich derivatives of single-valued mappings in $\mathbb{R}^2$.

For the purpose to find some formulas of Fréchet derivatives and Mordukhovich derivatives of single-valued mappings in Euclidean spaces, we provide examples of polynomial mappings, rational mappings, exponential mappings and logarithm mappings in $\mathbb{R}^2$. which includes. For each given mapping, we find it's Fréchet derivatives, Mordukhovich derivatives, and covering constants.

From the examples studied in this paper, we see the difficulty and complexity for finding the covering constants of single-valued mappings in Euclidean spaces. It is easy to understand that will be much more difficult and more complicated to find the covering constants of single-valued mappings in Banach spaces. If we consider Hilbert spaces to be special cases of Banach spaces, we have the following problems for consideration by interested readers.

(i) Find the principles for the existence of the Fréchet derivatives of single-valued mappings in Hilbert spaces;
(ii) Find the guidelines for calculating the Mordukhovich derivatives of single-valued mappings in Hilbert spaces;
(iii) Find an applicable algorithm to calculate the covering constants for single-valued mappings in Hilbert spaces;
(iv) Consider the above problems for set-valued mappings in Hilbert spaces;
(v) Consider the above problems for single valued and set-valued mappings in Asplund spaces.